\documentclass[12pt]{article}
\usepackage{amssymb, amsmath,epsfig,amsthm}
\begin{document}
 \righthyphenmin=2

\def \<{{\langle}}
\def\>{{\rangle}}
\def\Z{{\Bbb Z}}
\def\N{{\mathbb N}}
\def\NN{{\mu_1}}
\def\T{{\cal T}}
\def \D{{\cal D}}
\def \C{{\cal C}}
\def \B{{\cal B}}
\def \Di{{\cal D}}
\def \di{{\Bbb D}}
\def \og{{\Bbb O}}
\def \H{{\cal H}}
\def \I{{\cal I}}
\def \Tht{{\Theta}}
\def \F{{\cal F}}
\def \V{{\cal V}}
\def \a{{\alpha}}
\def \l{{\lambda}}
\def \t{{\beta}}
\def \bt{{\beta}}
\def \dt{{\delta}}
\def \De{{\Delta}}
\def \ga{{\gamma}}
\def \kp{{\varkappa}}
\def \sgn{{\operatorname{sgn}}}
\def \R{{\mathbb R}}
\def \CR{{\cal R}}
\def \A {{\cal A}}
\def \a{{\alpha}}
\def \m{{\mu}}
\def \M{{\cal M}}
\def \g{{\gamma}}
\def \bt{{\beta}}
\def \la{{\lambda}}
\def \La{{\Lambda}}
\def \fee{{\varphi}}
\def \e{{\epsilon}}
\def \s{{\sigma}}
\def \Sm{{\Sigma}}
\def \sm{{\sigma_m}}
\def \PLb{{P^\bot_L}}
\def \td{\tilde}
\newcommand{\wh}{\widehat}
\newcommand{\wt}{\widetilde}
\newcommand{\ol}{\overline}
\def \supp{{\operatorname{supp}}}
\def \sign{{\operatorname{sign}}}
\def \sp{{\operatorname{span}}}
\def \csp{{\overline{\operatorname{span}}}}
\def \Prj{{\operatorname{Proj}}}
\def \Dff{{\operatorname{Diff}}}
\def \Dst{{\operatorname{Dist}}}
\def\charf\#1{\chi_{\lower3pt\hbox{$\scriptstyle \#1$}}} 
\catcode`\@=11
\def\NoLogo{\let\logo@\empty}
\catcode`\@=\active \NoLogo

\newtheorem{Theor}{Theorem}
\newtheorem{Lemm}{Lemma}
\newtheorem{Corollar}{Collary}
\theoremstyle{remark}
\newtheorem{Remar}{Remark}

\begin{center}
{\Large On the optimality of Orthogonal Greedy Algorithm for
$M$-coherent dictionaries.}\footnote{This research is partially
supported by Russian Foundation for Basic Research
project 08-01-00799 and 09-01-12173} \\
 Eugene Livshitz
\end{center}
\begin{abstract}
We show that Orthogonal Greedy Algorithms (Orthogonal Matching
Pursuit) provides almost optimal approximation on the first
$[1/(20M)]$ steps for $M$-coherent dictionaries.
\end{abstract}

\section{Introduction.}

In this article we continue the research of convergence of greedy
algorithms with regards to dictionaries with small coherence (see
\cite{GMS}, \cite{GN}, \cite{Tr}, \cite{DET}, \cite{DET1},
\cite{TZh}, \cite{L}). The study of approximation by incoherent
dictionaries was mainly motivated by applications to compressed
sensing. In \cite{GMS}, \cite{Tr}, \cite{DET} it was shown that
Orthogonal Greedy Algorithm (Orthogonal Matching Pursuit) is
effective for signal recovering. In this article we discuss this
problem from the point of view of Approximation Theory.

Let us recall standard definition of Greedy Algorithms theory. We
say that a set $\D$ from a Hilbert space $H$ is a dictionary if
\begin{equation*}
 \phi\in\D\ \Rightarrow \|\phi\|=1,\mbox{ and }
\csp{\D}=H.
\end{equation*}
We study dictionaries with small values of coherence
\begin{equation}\label{cohM}
M:=\sup_{\phi,\psi\in\D,\ \phi\ne\psi}|\<\phi,\psi\>|.
\end{equation}
Dictionaries with coherence $M$ are called $M$-coherent.

{\sc Orthogonal Greedy Algorithm (OGA)} {\it Set  $f_0:=f\in H$,
$G_0^{OGA}(f,\Di) := 0$. For each $m\ge 0$ we inductively find
$g_{m+1}\in \D$ such that
\begin{equation*}
|\<f_m,g_{m+1}\>|= \sup_{g\in\D}|\<f_m,g\>|
\end{equation*}
and define}
\begin{equation*}
 G_{m+1}^{OGA}(f,\D):= \Prj_{\sp(g_{1},\ldots,g_{m+1})}(f),
\end{equation*}
\begin{equation*}
f_{m+1}:= f-G_{m+1}^{OGA}(f,\Di).
\end{equation*}

For a function $f\in H$ we define its best $m$-term approximation
\begin{equation*}
\sm(f):=\sm(f,\D):=\inf_{c_i\in\R,\phi_i\in\D, 1\le i\le m}
\|f-\sum_{i = 1}^mc_i\phi_i\|.
\end{equation*}

Following V.N.~Temlyakov we call inequalities connecting the error
of Greedy approximation and the best $m$-term approximation {\it
Lebesgue type inequalities}.

The first Lebesgue type inequality for Greedy Algorithms was
obtained by A.C.~Gilbert, M.~Muthukrishnan and J.~Strauss in
\cite{GMS}

{\bf Theorem A.} {\it For every $M$-coherent dictionary $\D$ and
any function $f\in H$ the inequality
\begin{equation*}
\|f-G^{OGA}_m(f,\D)\|=\|f_m\|\le 8m^{1/2}\sm(f)
\end{equation*}
holds for all $m$,
\begin{equation*}
1\le m\le \frac{1}{8\sqrt{2}M}-1.
\end{equation*}
}

This estimate was improved by J.~Tropp \cite{Tr} (see also
paper~\cite{DET} of D.\,L.~Donoho, M.~Elad and V.N.~Temlyakov.)

 {\bf Theorem B.}
{\it  For every $M$-coherent dictionary $\D$ and any function
$f\in H$
\begin{equation*}
\|f_m\|\le (1+6m)^{1/2}\sm(f),\mbox{ if }\ 1\le m\le \frac{1}{3M}.
\end{equation*}
}

D.\,L.~Donoho, M.~Elad and V.N.~Temlyakov \cite{DET1} dramatically
improved factor in front of $\s$.

{\bf Theorem C.} {\it  For every $M$-coherent dictionary $\D$ and
any function $f\in H$
\begin{equation*}
\|f_{\lfloor m\log m\rfloor}\|\le 24\sm(f),\mbox{ if }\ 1\le m\le
\frac{1}{20M^{2/3}}.
\end{equation*}
}

V.N.~Temlyakov and P.~Zheltov \cite{TZh} improved the upper border
for $m$ and proved two new Lebesgue type inequalities.

{\bf Theorem D.} {\it  For every $M$-coherent dictionary $\D$ and
any function $f\in H$
\begin{equation*}
\|f_{m\lfloor 2^{\sqrt{\log m}}\rfloor}\|\le 3\sm(f),\mbox{ if }\
m2^{\sqrt{2\log m}}\le \frac{1}{26M}.
\end{equation*}
}

{\bf Theorem E.} {\it  For every $M$-coherent dictionary $\D$, any
function $f\in H$ and any fixed $\dt>0$
\begin{equation*}
\|f_{m2^{\lceil\frac{1}{\dt}\rceil}}\|\le 3\sm(f),\mbox{ if }\
m\le
\left(\frac{1}{14M}\right)^{\frac{1}{1+\dt}}2^{-\lceil\frac{1}{\dt}\rceil}.
\end{equation*}
}

The aim of this article is to prove the following result.
\begin{Theor}\label{Th-cohOGA} For every $M$-coherent dictionary $\D$ and any
function $f\in H$ we have
\begin{equation*}
\|f-G^{OGA}_{2m}(f,\D)\|=\|f_{2m}\|\le 3\sm(f)
\end{equation*}
for all
\begin{equation*}
1\le m\le \frac{1}{20M}.
\end{equation*}
\end{Theor}

\section{Preliminary lemmas.}

By conditions of Theorem~\ref{Th-cohOGA} we have
\begin{equation}\label{MmM}
M\le mM \le 1/20.
\end{equation}

 We use several
standard lemmas.

\begin{Lemm}\label{Lm-max} For any $n$, $1\le n\le 2m$, and
\begin{equation*}
h=\sum_{i=1}^n c_i\phi_i,\  c_i\in\R, \phi_i\in\D,
\end{equation*}
we have
\begin{equation}\label{InnLeCoef}
\max_{1\le i\le n}|\<h,\phi_i\>|\le \max_{1\le i\le
n}|c_i|(1+2mM),
\end{equation}
\begin{equation}\label{InnGeCoef}
\max_{1\le i\le n}|\<h,\phi_i\>|\ge \max_{1\le i\le
n}|c_i|(1-2mM),
\end{equation}
\begin{equation}\label{CoefLeInn}
\max_{1\le i\le n}|c_i|\le \max_{1\le i\le n}|\<h,\phi_i\>|
(1+3mM),
\end{equation}
\end{Lemm}
\begin{proof}
Using (\ref{cohM}) we have for any $1\le i\le n$
\begin{multline*}
\<h,\phi_i\>=\<c_i\phi_i,\phi_i\> + \<\sum_{1\le j\le n,\ i\ne j}
c_j\phi_j,\phi_i\>\le\\
\le c_i+ (n-1)\left(\max_{1\le i\le n}|c_i|\right)M\le c_i +
\left(\max_{1\le i\le n}|c_i|\right)2mM.
\end{multline*}
Similarly
\begin{equation*}
\<h,\phi_i\>\ge c_i - \left(\max_{1\le i\le n}|c_i|\right)2mM.
\end{equation*}
The last two inequalities imply (\ref{InnLeCoef}) and
(\ref{InnGeCoef}). It follow form (\ref{MmM}) that
\begin{equation*}
(1-2mM)(1+3mM) = 1 + mM - 6(mM)^2 \ge 1+ mM(1-0.3)\ge 1.
\end{equation*}
To prove (\ref{CoefLeInn}) we estimate
\begin{equation*}
\max_{1\le i\le n}|c_i|\le \max_{1\le i\le n}|\<h,\phi_i\>|
(1-2mM)^{-1} \le \max_{1\le i\le n}|\<h,\phi_i\>| (1+3mM).
\end{equation*}
\end{proof}

As corollary we obtain
\begin{Lemm}\label{Lm-max1} Let $n\le 2m$, $h\in H$, $\phi_i\in\D$, $1\le i\le
n$. Suppose that
\begin{equation*}
\Prj_{\sp(\phi_1,\ldots,\phi_n)}(h)=\sum_{i = 1}^n c_i\phi_i.
\end{equation*}
Then \begin{equation*} \max_{1\le i\le n}|c_i|\le \max_{1\le i\le
n}|\<h,\phi_i\>| (1+3mM),
\end{equation*}
\end{Lemm}
\begin{proof}
Set
\begin{equation*}
h' = \Prj_{\sp(\phi_1,\ldots,\phi_n)}(h).
\end{equation*}
It's clear that
\begin{equation*}
\<h,\phi_i\> = \<h',\phi_i\>,\ 1\le i\le n.
\end{equation*}
Thus the lemma follows from inequality (\ref{CoefLeInn}) for $h'$.
\end{proof}

For $n\ge 1$ we define
\begin{equation}\label{dn-def}
d_n:=\<f_{n-1},g_n\>.
\end{equation}
Let numbers $x_{i,n}$, $n\ge 1$, $1\le i\le n$, satisfy the
equality
\begin{equation}\label{numbxin-def}
f_{n} = f_{n-1} - \sum_{i = 1}^n x_{i,n}g_i.
\end{equation}

\begin{Lemm}\label{Lm-xin} For any $n\le 2m$
we have
\begin{equation}\label{Ineqxin}
|x_{i,n}|\le M|d_n|(1+3mM),\ 1\le i\le n-1,
\end{equation}
\begin{equation}\label{Ineqxnn}
|x_{n,n}-d_n|\le M|d_n|(1+3mM).
\end{equation}
\end{Lemm}
\begin{proof}
By definition of OGA
\begin{equation}\label{inner-zero}
\<f_l,g_i\> = 0,\ 1\le i\le l,
\end{equation}
and
\begin{equation*}
f_{n}= f-G_{n}^{OGA}(f,\Di) = f -
\Prj_{\sp(g_{1},\ldots,g_{n})}(f).
\end{equation*}
Hence
\begin{equation}\label{fnrepr}
f_{n}= f_{n-1}- \Prj_{\sp(g_{1},\ldots,g_{n})}(f_{n-1}) = f_{n-1}
- d_ng_n - \Prj_{\sp(g_{1},\ldots,g_{n})}(f_{n-1} - d_ng_n).
\end{equation}
Using (\ref{cohM}) and (\ref{inner-zero}) we have for $h:=f_{n-1}
- d_ng_n$ and $1\le i\le n-1$
\begin{equation*}
|\<h,g_i\>| \le |\<f_{n-1},g_i\>| + |d_n\<g_i,g_n\>|\le M|d_n|,
\end{equation*}
\begin{equation*}
\<h,g_n\>=0.
\end{equation*}
Suppose that $x'_{i,n}$, $1\le i\le n$ satisfy
\begin{equation*}
\Prj_{\sp(g_{1},\ldots,g_{n})}(f_{n-1} - d_ng_n) =
\Prj_{\sp(g_{1},\ldots,g_{n})}(h) = \sum_{i=1}^n x'_{i,n}g_i.
\end{equation*}
By Lemma~\ref{Lm-max1}
\begin{equation}\label{xprinj}
|x'_{i,n}|\le M|d_n|(1+3mM),\ 1\le i\le n
\end{equation}
It follows from (\ref{fnrepr}) that
\begin{equation*}
 f_{n} = f_{n-1} -d_ng_n - \sum_{i=1}^n x'_{i,n}g_i =
 f_{n-1} - \sum_{i =
1}^n x_{i,n}g_i,
\end{equation*}
where $x_{i,n} = x'_{i,n}$, $1\le i\le n-1$ and $x_{n,n} =
d_n+x'_{n,n}$. This and (\ref{xprinj}) complete the proof.
\end{proof}

\begin{Lemm}\label{Lm-dnp1} For $1\le n\le 2m-1$ we have
\begin{equation*}
|d_{n+1}|\le |d_n|(1+1.25M).
\end{equation*}
\end{Lemm}
\begin{proof}
By definition of OGA
\begin{equation*}
|\<f_{n-1},g_{n+1}\>|\le |\<f_{n-1},g_{n}\>| = |d_{n}|,
\end{equation*}
Using Lemma~\ref{Lm-xin}, (\ref{cohM}) and (\ref{MmM})  we have
\begin{multline*}
|\<f_{n},g_{n+1}\>|\le |\<f_{n-1} - \sum_{i = 1}^n
x_{i,n}g_i,g_{n+1}\>|\le |\<f_{n-1},g_{n+1}\>| + \sum_{i=1}^n
|x_{i,n}\<g_i,g_{n+1}\>|\le\\
\le |d_n|+M\left(|x_{n,n}|+\sum_{i=1}^{n-1}|x_{i,n}|\right)
\le |d_n| +(n M|d_n|(1+3mM) + |d_n|) M \le\\
\le |d_n|\left(1 + (2mM(1+3mM) + 1)M\right) \le |d_n|(1+1.25M).
\end{multline*}
\end{proof}

\begin{Lemm}\label{Lm-dln} For any $1\le l\le n \le 2m$ we have
\begin{equation*}
|d_{n}|\le |d_l|\exp(2.5mM).
\end{equation*}
\end{Lemm}
\begin{proof}
Using Lemma~\ref{Lm-dnp1} we write
\begin{equation*}
|d_n| \le |d_l|(1+1.25M)^{n-l}\le
|d_l|\left(1+\frac{2.5mM}{2m}\right)^{2m}\le |d_l|\exp(2.5mM).
\end{equation*}
\end{proof}

\section{Notations.}
By the definition of the best $m$-term approximation there exist
$a_j\in\R$, $\psi_j\in\D$, $1\le j\le m$, and $v_0\in H$ such that
\begin{equation}\label{f-Represent}
f=f_0=\sum_{j=1}^m a_j\psi_j + v_0,\quad \<v_0,\psi_j\>=0,\ 1\le
j\le m,
\end{equation}
\begin{equation}\label{v0sm}
\|v_0\|\le 1.01\sm(f).
\end{equation}
Set
\begin{equation*}
L:=\sp(\psi_1,\ldots,\psi_m),\quad P_L(\cdot):=\Prj_L(\cdot),\quad
P^\bot_L(\cdot):=\Prj_{L^\bot}(\cdot),
\end{equation*}
\begin{equation*}
v_n:=P^\bot_L(f_n),\ 0\le n\le 2m.
\end{equation*}
Let numbers $a_{j,n}$ and $b_{j,n}$, $0\le n\le 2m$, $1\le j\le m$
satisfy equalities
\begin{equation}\label{fnai}
f_n = P_L(f_n) + P^\bot_L(f_n) = \sum_{j = 1}^m a_{j,n}\psi_j +
v_n.
\end{equation}
\begin{equation}\label{bin-def}
\sum_{j = 1}^m b_{j,n}\psi_j = P_L(f_0-f_n).
\end{equation}
Then
\begin{equation}\label{ainbin}
a_{j,n} = a_j - b_{j,n},\ 1\le j\le m,\ 1\le n\le 2m.
\end{equation}

 Define
\begin{equation*}
T_1:=\left\{i\in\{1,\ldots,2m\}:\
g_i\in\{\psi_j\}_{j=1}^{m}\right\}.
\end{equation*}
\begin{equation*}
T_2:=\{1,\ldots,2m\}\setminus T_1,
\end{equation*}
\begin{equation}\label{S1-def}
S_1:=\left\{j\in\{1,\ldots,m\}:\
\psi_j\in\{g_n\}_{n=1}^{2m}\right\},
\end{equation}
\begin{equation}\label{S2-def}
S_2:=\{1,\ldots,m\}\setminus S_1.
\end{equation}
For numbers $x_{i,n}$, $1\le i\le n\le 2m$, from
(\ref{numbxin-def}) and for $d_n$ from (\ref{dn-def}) we define
\begin{equation*}
x_n:=\sum_{1\le i\le n,\ i\in T_2} |x_{i,n}|,
\end{equation*}
\begin{equation*}
D:=\sum_{1\le n\le 2m,\ n\in T_2}d_n^2.
\end{equation*}

\section{Main lemmas.}

\begin{Lemm}\label{Lm-proj} Let $1\le i < n\le 2m$, $i,n\in T_2$.
Then we have
\begin{equation*}
|\<\PLb(g_n),g_i\>|\le 1.1M
\end{equation*}
\end{Lemm}
\begin{proof}
Let
\begin{equation*}
P_L(g_n)=\sum_{j= 1}^m c_j\psi_j.
\end{equation*}
Since $n\in T_2$ and
\begin{equation*}
g_n\ne\psi_j,\quad |\<g_n,\psi_j\>|\le M,\quad 1\le j\le m,
\end{equation*}
we get by Lemma~\ref{Lm-max1} that
\begin{equation*}
\max_{1\le j\le m}|c_j|\le M(1+3mM).
\end{equation*}
Therefore we have
\begin{multline*}
|\<\PLb(g_n),g_i\>| = |\<g_n-P_L(g_n),g_i\>|\le |\<g_n,g_i\>| +
|\<P_L(g_n),g_i\>|\le\\
\le M + |\<\sum_{j= 1}^m c_j\psi_j,g_i\>| \le M + m
\left(\max_{1\le j\le m}|c_j|\right) \max_{1\le j\le m}
|\<\psi_j,g_i\>|\le\\
\le
 M + (mM)M(1+3mM)\le 1.1 M.
\end{multline*}

\end{proof}

\begin{Lemm}\label{Lm-main1} Let $n\in T_1$ then
\begin{equation*}
x_n\le 0.1 D^{1/2}m^{-1/2},
\end{equation*}
\begin{equation*}
 \|v_n\|^2\le \|v_{n-1}\|^2 + 0.3DM.
\end{equation*}
\end{Lemm}
\begin{proof}
Let
\begin{equation*}
u_n:=\sharp\left(T_2\cap\{1,\ldots, n\}\right).
\end{equation*}
If $T_2\cap\{1,\ldots, n\} = \emptyset$ then $x_n=0$,
$v_n=v_{n-1}=v_0$ and nothing to prove, so we may assume that
$u_n\ge 1$. By Lemma~\ref{Lm-dln}
\begin{equation}\label{dn-ineq1}
|d_n|\le \exp(2.5mM)\min_{1\le i\le n,\ i\in T_2}|d_i|.
\end{equation}
On the other hand we have
\begin{equation*}
\left(\min_{1\le i\le n,\ i\in T_2}|d_i|\right)^2u_n\le \sum_{1\le
i\le n,\ i\in T_2}d_i^2 \le \sum_{1\le i\le 2m,\ i\in T_2}d_i^2 =
D.
\end{equation*}
Combining with (\ref{dn-ineq1}) we obtain
\begin{equation}\label{dnun}
|d_n|\le \exp(2.5mM)\left(\frac{D}{u_n}\right)^{1/2}.
\end{equation}
\begin{equation}\label{dnun1}
d_n^2u_n\le \exp(5mM)D.
\end{equation}
Applying Lemma~\ref{Lm-xin}, (\ref{MmM})  and (\ref{dnun}) we
write
\begin{multline}\label{x_n-desineq}
x_n=\sum_{1\le i\le n,\ i\in T_2} |x_{i,n}| = \sum_{1\le i\le n -
1,\ i\in T_2} |x_{i,n}|  \le M|d_n|(1+3mM) u_n\le\\
\le M(1+3mM)\exp(2.5mM)(Du_n)^{1/2} \le
M(1+3mM)\exp(2.5mM)(D2m)^{1/2}
=\\
= (2D)^{1/2}Mm^{1/2}(1+3mM)\exp(2.5mM)\le 0.1 D^{1/2}m^{-1/2}.
\end{multline}
We have that for any $l$, $1\le l\le m$ we have

\begin{equation*}
|\<\sum_{j = 1}^m a_{j,n-1}\psi_j,\psi_l\>|= |\<\sum_{j = 1}^m
a_{j,n-1}\psi_j +v_{n-1},\psi_l\>|\le |\<f_{n-1},\psi_l\>|\le
|d_n|.
\end{equation*}

 Then by Lemma~\ref{Lm-max} we get

\begin{equation}\label{ainm1}
\max_{1\le j\le m}|a_{j,n-1}|\le\left(\max_{1\le l\le m}|\<\sum_{j
= 1}^m a_{j,n-1}\psi_j,\psi_l\>|\right)(1+3mM)\le |d_n|(1+3mM).
\end{equation}
Define
\begin{equation}\label{defH}
h:= \sum_{1\le i\le n, i\in T_2}x_{i,n}g_i = \sum_{1\le i\le n-1,
i\in T_2}x_{i,n}g_i.
\end{equation}
According the definition of $v_n$ we have
\begin{equation*}
v_n=\PLb(f_n)=\PLb\left(f_{n-1}-\sum_{i=1}^n
x_{i,n}g_i\right)=v_{n-1}-\PLb(h),
\end{equation*}
\begin{multline}\label{vnsqest}
\|v_n\|^2=\|v_{n-1}-\PLb(h)\|^2 \le \|v_{n-1}\|^2 +
2|\<v_{n-1},\PLb(h)\>| +\|\PLb(h)\|^2 \le\\
\le \|v_{n-1}\|^2 + 2|\<v_{n-1},h\>| +\|h\|^2.
\end{multline}
By definition of OGA $\<f_{n-1},g_i\> = 0$, $1\le i\le n-1$,
therefore using (\ref{defH}) and (\ref{fnai})
\begin{equation*}
\<f_{n-1},h\> = 0,
\end{equation*}
\begin{multline*}
|\<v_{n-1},h\>| =|\<f_{n-1}-\sum_{j = 1}^m a_{j,n-1}\psi_j,h\>| =
\sum_{j = 1}^m |\<a_{j,n-1}\psi_i,h\>|\le\\
\le \sum_{j = 1}^m |a_{j,n-1}|\sum_{1\le i\le n-1,\ i\in
T_2}|\<\psi_j, x_{i,n}g_i\>|.
\end{multline*}
Applying (\ref{cohM}), (\ref{dnun1}), (\ref{ainm1}) and
Lemma~\ref{Lm-xin} we obtain
\begin{multline*}
|\<v_{n-1},h\>| \le \sum_{j = 1}^m |a_{j,n-1}|\sum_{1\le i\le
n-1,\
i\in T_2}|x_{i,n}\<\psi_j, g_i\>| \le\\
\le |d_n|(1+3mM) M|d_n|(1+3mM) u_n mM \le\\
\le (d_n^2u_n)(1+3mM)^2 mM^2 \le D\exp(5mM) (1+3mM)^2 mM^2.
\end{multline*}

\begin{multline*}
\|h\|^2 = \sum_{1\le i\le n,\ i\in T_2 } x_{i,n}^2\<g_i,g_i\> + 2
\sum_{1\le i,l\le n,\ i,l\in T_2,\ i\ne
l}x_{i,n}x_{l,n}\<g_i,g_l\>\le\\
\le M^2|d_n|^2(1+3mM)^2u_n +
2M^2|d_n|^2(1+3mM)^2u_n^2M \le\\
\le M^2|d_n|^2u_n(1+3mM)^2 + M^2|d_n|^2u_n(1+3mM)^24mM\le\\
\le DM^2\exp(5mM)(1+3mM)^2(1+4mM).
\end{multline*}
From (\ref{vnsqest}) and (\ref{MmM}) it follows that
\begin{multline*}
\|v_n\|^2 \le \|v_{n-1}\|^2 + 2|\<v_{n-1},h\>| +\|h\|^2 \le\\
\le \|v_{n-1}\|^2 + DM\exp(5mM)(1+3mM)^2(2mM+M+4mM^2)\le\\
\le \|v_{n-1}\|^2 + DM\exp(5mM)(1+3mM)^2(3mM+4(mM)^2) \le
\|v_{n-1}\|^2 + 0.3DM.
\end{multline*}

This estimate together with (\ref{x_n-desineq}) proves the lemma.
\end{proof}

\begin{Lemm}\label{Lm-main2} Let $n\in T_2$ then
\begin{equation*}
x_n\le 1.15|d_n|,
\end{equation*}
\begin{equation*}
 \|v_n\|^2\le \|v_{n-1}\|^2 - 0.6 d_n^2
\end{equation*}
\end{Lemm}
\begin{proof}
Applying Lemma~\ref{Lm-xin} we have
\begin{multline*}
x_n = \sum_{1\le i\le n,\ i\in T_2} |x_{i,n}| \le d_n + \sum_{1\le
i\le n,\ i\in T_2}M|d_n|(1+3mM)\le\\
\le |d_n|(1+2mM(1+3mM))\le |d_n|(1+3mM)\le 1.15|d_n|.
\end{multline*}
By Lemma~\ref{Lm-max} we have
\begin{multline}\label{ajnm2}
\max_{1\le j \le m}|a_{j,n-1}|\le(1+3mM)\max_{1\le l\le m}
\left|\<\sum_{j=1}^m a_{j,n-1}\psi_j,\psi_l\>\right|=\\
= (1+3mM)\max_{1\le l\le m} \left|\<\sum_{j=1}^m
a_{j,n-1}\psi_j+v_{n-1},\psi_l\>\right|= \\
=(1+3mM)\max_{1\le l\le m}|\<f_{n-1},\psi_l\>| \le (1+3mM)|d_{n}|.
\end{multline}
Therefore
\begin{multline}\label{vnm1gn}
|\<v_{n-1},g_n\> - d_n| = \left|\<f_{n-1}-\sum_{j=1}^m
a_{j,n-1}\psi_j,g_n\>-d_n\right|=\\
= \left|\<f_{n-1},g_n\> - \sum_{j=1}^m
a_{j,n-1}\<\psi_j,g_n\>-d_n\right| \le \left| \sum_{j=1}^m
a_{j,n-1}\<\psi_j,g_n\>\right|\le\\
\le \left(\max_{1\le j \le m}|a_{j,n-1}|\right)m\max_{1\le j\le m
}|\<\psi_j,g_n\>|\le
 (1+3mM)|d_{n}|mM.
\end{multline}
 Set
\begin{equation*}
v_n':=\PLb(f_{n-1}-x_{n,n}g_n).
\end{equation*}
Using Lemma~\ref{Lm-xin}, (\ref{MmM}) and (\ref{vnm1gn}) we
estimate
\begin{multline*}
2x_{n,n}\<v_{n-1},g_n\> =2 (d_{n} + (x_{n,n} - d_n))
(d_n+(\<v_{n-1},g_n\>-d_n))\ge\\
\ge 2 (|d_{n}| - M(1+3mM)|d_{n}|) (|d_{n}|- (1+3mM)|d_{n}|mM)
\ge\\
 \ge 2d_n^2-2|d_n|^2(1+3mM)(M+mM)\ge 2d_n^2-4|d_n|^2(1+3mM)mM,
\end{multline*}
\begin{multline}\label{vprsqnorm}
\|v_n'\|^2\le\|\PLb(f_{n-1}-x_{n,n}g_n)\|^2=
\|v_{n-1}-x_{n,n}\PLb(g_n)\|^2\le\\
\le \|v_{n-1}\|^2
 -2x_{n,n}\<v_{n-1},\PLb(g_n)\>+x_{n,n}^2\|\PLb(g_n)\|^2\le\\
 \le
 \|v_{n-1}\|^2-2x_{n,n}\<v_{n-1},g_n\>+x_{n,n}^2\le\\
 \le
 \|v_{n-1}\|^2-2d_n^2+4|d_n|^2(1+3mM)mM +
 (|d_n|+M|d_n|(1+3mM))^2\le\\
 \le
 \|v_{n-1}\|^2 - 0.65d_n^2.
\end{multline}
Similar to the proof of Lemma~\ref{Lm-main1} we define
\begin{equation*}
h:= \sum_{1\le i\le n-1, i\in T_2}x_{i,n}g_i.
\end{equation*}
Equalities $\<f_{n-1},g_i\> = 0$, $1\le i\le n-1$ imply that
\begin{equation*}
\<f_{n-1},h\> = 0.
\end{equation*}
Using (\ref{fnai})  we have
\begin{multline*}
|\<v'_n,h\>| = |\<\PLb(f_{n-1})-x_{n,n}\PLb(g_n),h\>|=
|\<v_{n-1}-x_{n,n}\PLb(g_n),h\>| = \\
=|\<f_{n-1}-\sum_{j = 1}^m a_{j,n-1}\psi_j -x_{n,n}\PLb(g_n),h\>|
\le \sum_{j = 1}^m |\<a_{j,n-1}\psi_i,h\>|
+|x_{n,n}\<\PLb(g_n),h\>| \le\\
\le \sum_{j = 1}^m |a_{j,n-1}|\sum_{1\le i\le n-1,\ i\in
T_2}|\<\psi_j, x_{i,n}g_i\>|+\sum_{1\le i\le n-1,\ i\in
T_2}|x_{n,n}x_{i,n}\<\PLb(g_n),g_i\>|.
\end{multline*}
Applying (\ref{cohM}), (\ref{ajnm2}), Lemma~\ref{Lm-xin} and
Lemma~\ref{Lm-proj} we continue
\begin{multline}\label{vprnh}
|\<v'_n,h\>| \le \sum_{j = 1}^m |a_{j,n-1}|\sum_{1\le i\le n-1,\
i\in T_2}|x_{i,n}||\<\psi_j, g_i\>|+\sum_{1\le i\le n-1,\ i\in
T_2}|x_{n,n}x_{i,n}\<\PLb(g_n),g_i\>| \le\\
\le \max_{1\le j \le m}|a_{j,n-1}|\max_{1\le i\le n-1,\ i\in T_2}
|x_{i,n}|\sum_{j = 1}^m\sum_{1\le i\le n-1,\ i\in T_2}
M+\\
+|d_n|(1+M(1+3mM))\max_{1\le i\le n-1,\ i\in T_2}
|x_{i,n}|\sum_{1\le
i\le n-1,\ i\in T_2}1.1M\le\\
\le d_n^2(1+3mM)^2 nmM^2  + 1.1d_n^2(1+M(1+3mM))(1+3mM)nM^2\le\\
\le d_n^2(1.15)^2 2(mM)^2  + 1.1d_n^2(1+mM(1.15))(1.15)2mM^2 \le
0.014 d_n^2,
\end{multline}
\begin{multline}\label{hsqnorm}
\|h\|^2 = \sum_{1\le i\le n-1,\ i\in T_2 } x_{i,n}^2\<g_i,g_i\> +
2 \sum_{1\le i,l\le n-1,\ i,l\in T_2,\ i\ne
l}x_{i,n}x_{l,n}\<g_i,g_l\>\le\\
\le M^2|d_n|^2(1+3mM)^2n + 2M^2|d_n|^2(1+3mM)^2n^2M \le\\
\le d_n^2(1+3mM)^2M(2Mm+8(Mm)^2)\le 0.008d_n^2.
\end{multline}

Using  definitions of $v_n'$ and $h$ we write
\begin{multline*}
v_n=\PLb(f_n)=\PLb\left(f_{n-1}- \sum_{i=1}^{n} x_{i,n}g_i\right)=
\PLb(f_{n-1}-x_{n,n}g_n)-\PLb\left(\sum_{i=1}^{n-1}
x_{i,n}g_i\right)=\\
=v'_n-\PLb\left(\sum_{1\le i\le n-1,\ i\in T_2}
x_{i,n}g_i\right)=v'_n-\PLb(h),
\end{multline*}
\begin{equation*}
\|v_n\|^2 = \|v'_n\|^2-2\<v'_n,\PLb(h)\> + \|\PLb(h)\|^2\le
\|v'_n\|^2+2|\<v'_n,h\>|+\|h\|^2.
\end{equation*}

Applying (\ref{vprsqnorm}), (\ref{vprnh}) and (\ref{hsqnorm}) we
obtain
\begin{multline*}
\|v_n\|^2 \le \|v'_n\|^2+2|\<v'_n,\PLb(h)\>|+\|h\|^2\le\\
\le \|v_{n-1}\|^2 - 0.65 d_n^2 + 2(0.014 d_n^2)+ 0.008d_n^2 \le
\|v_{n-1}\|^2 - 0.6d_n^2
\end{multline*}

\end{proof}

\begin{Lemm}\label{Lm-main3}
We have
\begin{equation*}
\sum_{n=1}^{2m} x_n\le 2D^{1/2}m^{1/2}.
\end{equation*}
\end{Lemm}
\begin{proof}
Using Cauchy inequality, Lemma~\ref{Lm-main2} and
Lemma~\ref{Lm-main1} we get
\begin{multline*}
\sum_{1\le n \le 2m,\ n\in T_2} x_n\le \sum_{1\le n\le 2m,\ n\in
T_2}
1.15|d_n|\le \\
\le 1.15 \left(\sum_{1\le n\le 2m,\ n\in
T_2}d_n^2\right)^{1/2}(2m)^{1/2}\le 1.7D^{1/2} m^{1/2},
\end{multline*}
\begin{multline*}
\sum_{n=1}^{2m} x_l = \sum_{1\le n\le 2m,\ n\in T_1} x_n +
\sum_{1\le n\le 2m,\ n\in T_2} x_n\le \sharp T_1\max_{1\le n\le
2m,\ n\in T_1} x_n+\sum_{1\le n\le 2m,\ n\in T_2}x_n\le\\
\le m (0.1D^{1/2}m^{-1/2}) + 1.7D^{1/2} m^{1/2}\le
2D^{1/2}m^{1/2}.
\end{multline*}
\end{proof}

\begin{Lemm}\label{Lm-main4} We have
\begin{equation*}
D^{1/2}\le 1.33\sm(f),
\end{equation*}
\begin{equation*}
\|v_{2m}\|\le \|v_0\|.
\end{equation*}

\end{Lemm}
\begin{proof}
Applying (\ref{v0sm}), Lemma~\ref{Lm-main1} and
Lemma~\ref{Lm-main2} we write
\begin{multline*}
(1.01\sm(f))^2\ge \|v_0\|^2\ge \|v_0\|^2-\|v_{2m}\|^2 =
\sum_{n=1}^{2m}
(\|v_{n-1}\|^2-\|v_n\|^2) =\\
= \sum_{1\le n\le 2m,\ n\in T_1} (\|v_{n-1}\|^2-\|v_n\|^2) +
\sum_{1\le n\le 2m,\ n\in T_2} (\|v_{n-1}\|^2-\|v_n\|^2)\ge\\
\ge \sharp T_1\left(-0.3DM\right)+\sum_{1\le n\le 2m,\ n\in
T_2}0.6d_n^2 \ge m(-0.3DM)+0.6D\ge 0.58D .
\end{multline*}
Hence
\begin{equation*}
D^{1/2}\le 1.01(0.58)^{-1/2}\sm(f)\le 1.33\sm(f).
\end{equation*}
\end{proof}

In the next lemma we use definitions (\ref{bin-def}),
(\ref{S1-def}) and (\ref{S2-def}).

\begin{Lemm}\label{Lm-main5} For any $1\le n\le 2m$ and $j\in S_2$
\begin{equation*}
|b_{j,n}| \le  0.12 D^{1/2}m^{-1/2}.
\end{equation*}
\end{Lemm}
\begin{proof}
By definition (\ref{numbxin-def})
\begin{multline*}
f_n-f_0=\sum_{l=1}^n\sum_{i=1}^l
x_{i,l}g_i=\sum_{i=1}^ng_i\left(\sum_{l=i}^n x_{i,l}\right) =\\
= \sum_{1\le i\le n,\ i\in T_1}g_i\left(\sum_{l=i}^n
x_{i,l}\right)+ \sum_{1\le i\le n,\ i\in T_2}g_i\left(\sum_{l=i}^n
x_{i,l}\right).
\end{multline*}
Assume that numbers $\wh b_{j,n}$ and $\wt b_{j,n}$, $1\le j\le
m$, $1\le n\le 2m$ satisfy
\begin{equation*}
\sum_{j=1}^m \wh b_{j,n}\psi_j=P_L\left(\sum_{1\le i\le n,\ i\in
T_1}g_i\left(\sum_{l=i}^n x_{i,l}\right)\right),
\end{equation*}
\begin{equation*}
\sum_{j=1}^m \wt b_{j,n}\psi_j=P_L\left(\sum_{1\le i\le n,\ i\in
T_2}g_i\left(\sum_{l=i}^n x_{i,l}\right)\right).
\end{equation*}
It follows from (\ref{bin-def}) that
\begin{equation}\label{bjn-eq}
b_{j,n} = \wh b_{j,n} +\wt b_{j,n},\ 1\le j\le m,\ 1\le n\le 2m.
\end{equation}
It's clear that
\begin{equation*}
P_L\left(\sum_{1\le i\le n,\ i\in T_1}g_i\left(\sum_{l=i}^n
x_{i,l}\right)\right) = \sum_{1\le i\le n,\ i\in
T_1}g_i\left(\sum_{l=i}^n x_{i,l}\right)
\end{equation*}
and therefore
\begin{equation}\label{whbjn}
\wh b_{j,n} = 0,\ j\in S_2,\ 1\le n\le 2m.
\end{equation}
Set
\begin{equation*}
h = \sum_{1\le i\le n,\ i\in T_2}g_i\left(\sum_{l=i}^n
x_{i,l}\right).
\end{equation*}
By Lemma~\ref{Lm-main3} we estimate for each $j$, $1\le j\le m$,
\begin{multline*}
|\<h,\psi_j\>|\le\sum_{1\le i\le n,\ i\in
T_2}\<g_i,\psi_j\>\left(\sum_{l=i}^n |x_{i,l}|\right)\le\\
\le M\sum_{l=1}^n\sum_{1\le i\le l,\ i\in T_2}|x_{i,l}|\le
M\sum_{l=1}^n x_l  \le M\sum_{l=1}^{2m} x_l \le 2D^{1/2}m^{1/2}M.
\end{multline*}
According Lemma~\ref{Lm-max1} we have for $1\le j\le m$ and $1\le
n\le 2m$
\begin{equation}\label{wtbjn}
|\wt b_{j,n}| \le 2D^{1/2}m^{1/2}M(1+3mM).
\end{equation}
Combining (\ref{MmM}), (\ref{bjn-eq}), (\ref{whbjn}) and
(\ref{wtbjn}) we obtain
\begin{equation*}
|b_{j,n}| \le  2D^{1/2}m^{1/2}M(1+3mM)\le  0.12 D^{1/2}m^{-1/2} ,\
j\in S_2,\ 1\le n\le 2m.
\end{equation*}
\end{proof}

\section{Proof of Theorem \ref{Th-cohOGA}.}

First we estimate $|a_{j,2m}|$, $j\in S_2$. For each $1\le n\le
2m$ by Lemma~\ref{Lm-max}  and Lemma~\ref{Lm-main5} we have
\begin{multline*}
|d_n| \ge \max_{1\le j\le m}|\<f_{n-1},\psi_j\>| = \max_{1\le j\le
m}|\<P_L(f_{n-1}),\psi_j\>| \ge (1-2mM)\max_{1\le j\le m}
|a_{j,n-1}|\ge\\
\ge(1-2mM)\max_{j\in S_2} |a_{j,n-1}| \ge (1-2mM) \max_{j\in S_2}
|a_j - b_{j,n-1}|\ge\\
\ge (1-2mM)\left(\max_{j\in S_2}|a_j| - 0.12
D^{1/2}m^{-1/2}\right).
\end{multline*}
Since $\sharp T_1\le m$ and $\sharp T_2\ge m$ we get
\begin{equation*}
D=\sum_{1\le i\le 2m,\ i\in T_2}d_n^2\ge
m\left((1-2mM)\left(\max_{j\in S_2}|a_j| - 0.12
D^{1/2}m^{-1/2}\right)\right)^2.
\end{equation*}
Hence
\begin{equation*}
\left(\max_{j\in S_2}|a_j| - 0.12 D^{1/2}m^{-1/2}\right)\le
D^{1/2}m^{-1/2}(1+3mM),
\end{equation*}
\begin{equation*}
\max_{j\in S_2}|a_j| \le 1.15D^{1/2}m^{-1/2} + 0.12
D^{1/2}m^{-1/2}= 1.27D^{1/2}m^{-1/2}.
\end{equation*}
Then by (\ref{ainbin}) and Lemma~\ref{Lm-main5} for any $j\in S_2$
we obtain
\begin{equation*}
|a_{j,2m}| = |a_j-b_{j,2m}|\le 1.27D^{1/2}m^{-1/2} + 0.12
D^{1/2}m^{-1/2}\le 1.4D^{1/2}m^{-1/2}.
\end{equation*}
We use well known inequality (see, for example, Lemma 2.1 from
\cite{DET1})
\begin{equation*}
\left\|\sum_{j\in S_2}a_{j,2m}\psi_j\right\|^2\le\left(\sum_{j\in
S_2}a_{j,2m}^2\right)(1+mM)\le
m\left(1.4D^{1/2}m^{-1/2}\right)^2(1.05) \le 2.06D.
\end{equation*}
Using the definition of OGA, (\ref{fnai}), Lemma~\ref{Lm-main4}
and (\ref{v0sm}) we estimate
\begin{multline*}
\|f_{2m}\|= \min_{c_i,\ 1\le i\le 2m}
\|f_{2m}-\sum_{i=1}^{2m}c_ig_i\|=
 \min_{c_i,\ 1\le i\le
2m}\|\sum_{j=1}^ma_{j,2m}\psi_j+v_{2m}-\sum_{i=1}^{2m}c_ig_i\|\le\\
\le \min_{c_i,\ 1\le i\le 2m}
\|\sum_{j=1}^ma_{j,2m}\psi_j-\sum_{i=1}^{2m}c_ig_i\|
+\|v_{2m}\|\le\\
\le \min_{c_l,\ l\in S_1} \|\sum_{j=1}^ma_{j,2m}\psi_j-\sum_{l\in
S_1}c_l\psi_l\| +\|v_{0}\|
\le \left\|\sum_{j\in S_2}a_{j,2m}\psi_j\right\| + \|v_0\| \le\\
\le (2.06D)^{1/2} + 1.01\sm(f)\le (2.06)^{1/2} 1.33\sm(f) +
1.01\sm(f)\le 3 \sm(f).
\end{multline*}
This completes the proof. $\square$

The author is grateful professor V.N.~Temlyakov and professor
S.V.~Konyagin for useful discussions.

\end{document}